\documentclass[11pt]{amsart}
\usepackage{mathptmx, amssymb, amsthm, amsmath, mathtools, enumerate, colonequals}

\usepackage{xcolor, color}
\usepackage[margin=1.25in]{geometry}
\definecolor{chianti}{rgb}{0.6,0,0}
\definecolor{meretale}{rgb}{0,0,.6}
\definecolor{leaf}{rgb}{0,.35,0}
\usepackage[colorlinks=true, pagebackref, hyperindex, citecolor=meretale, urlcolor=leaf, linkcolor=chianti]{hyperref}

\usepackage{leftindex}
\usepackage{graphicx,float}

\input{xy}
\xyoption{all}
\usepackage{tikz}

\usepackage{comment}
\usepackage{setspace}

\numberwithin{equation}{section}

\newtheorem{theorem}{Theorem}[section]
\newtheorem*{theoremA}{Theorem A}
\newtheorem*{theoremB}{Theorem B}
\newtheorem{lemma}[theorem]{Lemma}
\newtheorem{prop}[theorem]{Proposition}
\newtheorem{corollary}[theorem]{Corollary}

\theoremstyle{definition}
\newtheorem{defn}[theorem]{Definition}
\newtheorem{example}[theorem]{Example}

\theoremstyle{remark}
\newtheorem{remark}[theorem]{Remark}
\AtEndEnvironment{remark}{\null\hfill\ensuremath{\diamond}}

\numberwithin{equation}{subsection}

\theoremstyle{definition}

\usepackage [english]{babel}
\usepackage [autostyle, english = american]{csquotes}
\MakeOuterQuote{"}

\DeclareMathOperator{\Hom}{Hom}
\DeclareMathOperator{\RHom}{RHom}

\DeclareMathOperator{\Ext}{Ext}

\DeclareMathOperator{\Spec}{Spec}
\DeclareMathOperator{\Supp}{Supp}

\DeclareMathOperator{\la}{\langle}
\DeclareMathOperator{\ra}{\rangle}
\DeclareMathOperator{\injdim}{inj.dim}

\newcommand{\frob}{\operatorname{Frob}}
\newcommand{\ufrob}{\frob^{\operatorname{unit}}}
\newcommand{\cart}{\operatorname{Cart}}
\newcommand{\ucart}{\cart^{\operatorname{unit}}}

\newcommand{\cryscart}{\cart^{\operatorname{crys}}}
\newcommand{\unit}{\operatorname{unit}}

\renewcommand{\to}[1][]{\xrightarrow{\ #1\ }}

\renewcommand{\phi}{\varphi}

\renewcommand{\theta}{\vartheta}

\renewcommand{\epsilon}{\varepsilon}

\makeatletter
\@namedef{subjclassname@2020}{%
  \textup{2020} Mathematics Subject Classification}
\makeatother

\title{On the injective dimension of unit Cartier and Frobenius modules}
\author{Manuel Blickle, Daniel Fink, Alexandria Wheeler and Wenliang Zhang}

\address{Institut für Mathematik, Johannes Gutenberg-Universität Mainz, 55126 Mainz, Germany}
\email{blicklem@uni-mainz.de, dfink@uni-mainz.de} 
\thanks{M.B. and D.F. are partially supported by the Deutsche Forschungsgemeinschaft (DFG, German Research Foundation) – Project-ID 444845124 – TRR 326.}
\address{Department of Mathematics, Statistics, and Computer Science, University of Illinois at Chicago, Chicago, IL 60607}
\email{awheel23@uic.edu, wlzhang@uic.edu}
\thanks{A.W. and W.Z. are partially supported by NSF through DMS-1752081.}

\begin{document}

\begin{abstract}
Let $R$ be a regular $F$-finite ring of prime characteristic $p$. We prove that the injective dimension of every unit Frobenius module $M$ in the category of unit Frobenius modules is at most $\dim(\Supp_R(M))+1$. We further show that for unit Cartier modules the same bound holds over any noetherian $F$-finite ring $A$ of prime characteristic $p$. This shows that $\dim A+1$ is a uniform upper bound for the injective dimension of any unit Cartier module over a noetherian $F$-finite ring $A$.
\end{abstract}

\maketitle

\section{Introduction}

In the seminal paper \cite{LyubeznikFModulesApplicationsToLocalCohomology}, Lyubeznik introduced the theory of \emph{unit Frobenius modules}\footnote{in op. cit. these are called \emph{$F$-modules}, but we follow the notation of \cite{EmersonKisin,BhattLurie}} which was inspired by earlier work including \cite{PeskineSzpiro}, \cite{HartshorneSpeiserLocalCohomologyInCharacteristicP} and \cite{HunekeSharp}. By definition, a \emph{Frobenius module} over an $\mathbb{F}_p$-algebra $R$ is an $R$-module equipped with a structural $R$-module homomorphism $M \to F_*M$, which represents a left-action of the Frobenius $F$ on $M$. The Frobenius module $M$ is called \emph{unit Frobenius module} if the map $F^*M\to M$ associated with the structural map is an isomorphism. The theory of unit Frobenius modules has found numerous applications in the study of rings/varieties in prime characteristic $p$ and has generated a fruitful line of research ({\it cf.} \cite{BlickleDModonRFMod, EmersonKisin, MontanerBlickleLyubeznik2005, LyubeznikSomeLCInvariants, LyubeznikVanishingLCCharp, HochsterFModules, Blickle.MinimalGamma, BBCartierCrstals, BBLSZ1, ZhangLyubeznikNumbersProjectiveSchemes, MaGlobalDimension, MaZhangEulerian, LyubeznikSinghWalther, ZhangInjDimFandD, KLSZ_D_F_length, KatzmanZhangSupport, SwitalaZhangDichotomy,BhattLurie, Zhang_F_e_module}). 

In \S\ref{some generalities} we will recall the basic setup of the theory. One particular result in \cite{LyubeznikFModulesApplicationsToLocalCohomology} is the following theorem.
\begin{theorem}[Lyubeznik]
\label{injdim as R-mod}
Let $R$ be a noetherian regular ring of prime characteristic $p$. For every unit Frobenius module $M$ one has
\[\injdim_R (M) \leq \dim(\Supp_R (M))\]
i.e. the injective dimension of the \emph{underlying $R$-module} of $M$ is bounded by the dimension of its support. 
 \end{theorem}
 
In \cite{HochsterFModules}, Hochster showed that the category of unit Frobenius modules itself has enough injective objects. Subsequently Ma (\cite{MaGlobalDimension}) showed that, assuming $R$ is a regular $F$-finite ring, the injective dimension of a unit Frobenius module $M$ over $R$ in the category $\ufrob_R$ of unit Frobenius modules itself is bounded by $\dim R+1$:
 \begin{theorem}[Ma]
 \label{injdim as F-mod}
Let $R$ be an $F$-finite noetherian regular ring of prime characteristic $p$ that admits a canonical module $\omega_R$ such that $\omega_R\cong F^!\omega_R$. Then 
    \[\injdim_{\ufrob_R}(M)\leq \dim (R)+1\]
for every unit Frobenius module $M$.
\end{theorem}
 

The first main result obtained in this paper is a strengthening of Theorem \ref{injdim as F-mod} in the spirit of Theorem~\ref{injdim as R-mod}.
\begin{theoremA}[=Theorem \ref{injdim of F-module}]
\label{Theorem A}
Let $R$ be a regular $F$-finite ring of prime characteristic $p$. Then 
\[\injdim_{\ufrob_R} (M) \le \operatorname{inj.dim}_R(M)+1 \leq \dim(\Supp_R (M))+1\]
for every unit Frobenius module $M$.
\end{theoremA}
As shown in \cite{MaGlobalDimension} tensoring with the canonical module $\omega_R$ yields an equivalence between the category of unit Frobenius modules over $R$ and the category of \emph{unit Cartier modules} over $R$. Similar to Frobenius modules, a \emph{Cartier module} is an $R$-module $N$ together with a structural map $N \to F^{\flat}N$, encoding a right action of the Frobenius $F$. A Cartier module $N$ is called \emph{unit Cartier module} if its structural map is an isomorphism ({\it cf.} \S\ref{some generalities}).  

Recently the two first named authors gave in \cite[Theorem~1]{BlickleFink} an extension of Ma's Theorem~\ref{injdim as F-mod} for unit Cartier modules over any $F$-finite (not necessarily regular) ring:
\begin{theorem}
\label{thm: BF}
Let $A$ be an noetherian $F$-finite ring of prime characteristic $p$ and let $R\to A$ be a surjection from an $F$-finite regular ring of prime characteristic $p$. Then 
\[\injdim_{\ucart_A} (N)\leq \dim (R) +1\]
for every unit Cartier module $N$ over $A$. 
\end{theorem}
The shortcoming of this result from op. cit. was that the obtained bound is not intrinsic to $A$; it depends on a chosen surjection from a regular ring $R \to A$. We remove the dependence on this choice in our second main result:
\begin{theoremB}[=Theorem \ref{thm: global dim of Cartier crystals2}]
Let $A$ be a noetherian $F$-finite ring of prime characteristic $p$. Then 
\[\injdim_{\ucart_A} (N)\leq \dim(\Supp_A (N))+1\]
for every unit Cartier module $N$. In particular, the global injective dimension of the category of unit Cartier modules is at most $\dim A+1$.
\end{theoremB}
Furthermore, the bound in Theorem B is optimal ({\it cf.} Remark \ref{bound is sharp for CM ring}).

The article is organized as follows: in \S\ref{some generalities} we recall the necessary notions and results on unit Cartier modules and unit Frobenius modules. Then we show in \S\ref{inequalities inj dimension} the following inequality and equality:
\begin{itemize}
\item For any (unit) Cartier module $N$ over $A$, an arbitrary ring of positive characteristic $p$, we have: 
    \[ \injdim_{\cart_A} (N) \leq \injdim_{A} (N)+1  \]
\item For any unit Cartier module $N$  over a regular $F$-finite ring $R$ of positive characteristic $p$, we have:
    \[ \injdim_{\ucart_R} (N) = \injdim_{\cart_R} (N) \]
\end{itemize}
Combining these, we conclude the (equivalent) unit Cartier module version of Theorem A. In the final chapter \S\ref{KashiwaraForUnit} we adapt the methods from (c) in \cite{BlickleFink} to the setting of unit Cartier modules to derive the proof of Theorem B.


\section{Preliminaries on Frobenius and Cartier modules}
\label{some generalities}

In this section, we collect some basic background material on the functor $(-)^!$, unit Cartier modules and unit Frobenius modules. A more thorough treatment can be found in \cite{LyubeznikFModulesApplicationsToLocalCohomology, MaGlobalDimension, BBCartierCrstals, BlickleFink}, among others. 

Each noetherian ring $A$ of prime characteristic $p$ is equipped with a ring endomorphism $F: A\xrightarrow{a\mapsto a^p}A$ called the Frobenius endomorphism. To distinguish the source and target rings of the Frobenius, we will denote the target ring by $F_*A$ and its element by $F_*a$. Note that $F_*A$ is an $A$-algebra with the same ring structure as $A$ and its $A$-algebra structure is induced by the Frobenius: $b\cdot (F_*a)=F_*(b^pa)$ for every $b\in A$ and $F_*a\in F_*A$.
\subsection*{The functor $(-)^!$} 


Given any homomorphism $\varphi:A\to B$ of commutative rings, the restriction of scalars functor $\varphi_*$ has a left adjoint $\varphi^*$ and also a right adjoint $\varphi^{\flat}$ given by
\[
\varphi^{\flat}(M)=\Hom_A(B,M)\ {\rm and\ }\varphi^{\flat}(M\xrightarrow{f}N)=\Hom_A(B,M)\xrightarrow{g\mapsto g\circ f}\Hom_A(B,N)
\]
where the $B$-module structure on $\varphi^{\flat} (M)=\Hom_A(B,M)$ is given by pre-multiplication. The adjunction between $\varphi_*$ and $\varphi^{\flat}$ is an isomorphism of $A$-modules 
\[
    \Hom_A(\varphi_*M,N) \cong \varphi_*\Hom_B(M,\varphi^\flat N)
\]
which follows from $\otimes-\Hom$ adjunction. Note that $(-)^{\flat}$ is a 2-functor, i.e. it is compatible with composition of ring homomorphisms. As a right adjoint, $(-)^{\flat}$ is left exact, and we denote its derived functor by $\varphi^!(-)=\RHom_A(B,-)$.  


\begin{remark}
In \cite[Remark~13.6]{GabberTStructure}, Gabber showed that every $F$-finite noetherian ring $A$ of prime characteristic $p$ is a homomorphic image of a noetherian regular ring $R$ equipped with a finite $p$-basis. Such a regular ring with $p$-basis has a canonical choice of a canonical module, namely the top exterior power of the module of Kähler differentials $\omega_R = \bigwedge^{n} \Omega_R$, where $n$ is the rank of $\Omega_R$. It can easily be verified that the classical Cartier operator induces a canonical isomorphism $\kappa: \omega_R \to F^!\omega_R$. Shifting, we obtain a dualizing complex $\omega_R^\bullet = \omega_R[n]$ together with a canonical isomorphism $\omega^\bullet_R \to F^!\omega_R^\bullet$. If $\psi: R \to A$ denotes the quotient map, defining $\omega^\bullet_A=\psi^!\omega^\bullet_R$ yields via the 2-functoriality of $(-)^!$ the quasi-isomorphism
\[
    \omega^\bullet_A = \psi^!\omega^\bullet_R \to[\psi^!\kappa] \psi^!F_R^!\omega^\bullet_R \simeq (F_R \circ \psi)^!\omega^\bullet_R \simeq (\psi \circ F_A)^!\omega^\bullet_R \simeq F_A^!\psi^!\omega^\bullet_R \simeq F_A^!\omega^\bullet_A.
\]
This means that any noetherian $F$-finite ring $A$ has a dualizing complex $\omega^\bullet_A$ which is equipped with a quasi-isomorphism $\omega_A^\bullet \to[\simeq] F^!\omega^\bullet_A$ 

Furthermore, if $A$ is Cohen-Macaulay then the dualizing complex $\omega_A^\bullet$ is quasi-isomorphic to the $A$-module $\omega_A:=\Ext^c_R(A,\omega_R)$ (where $c:=\dim(R)-\dim(A)$) placed in appropriate homological degree. Consequently, we get an isomorphism 
\[\omega_A \to[\cong] F^\flat_A\omega_A \]
for the canonical module $\omega_A$ of $A$. Hence $\omega_A$ is a unit Cartier module.

In particular, the hypotheses in Ma's Theorem \ref{injdim as F-mod} is satisfied for all noetherian $F$-finite (regular) rings. In an upcoming paper \cite{BBST} it is shown that this does not depend on the presentation of $A$ as the quotient of a regular ring and also works globally for all noetherian $F$-finite schemes.
\end{remark}

\subsection*{(Unit) Cartier and (unit) Frobenius modules}

For a commutative ring $A$ of prime characteristic $p$, the category of Frobenius modules $\frob_A$ and the category of Cartier modules $\cart_A$ over $A$, are simply the categories of left and right modules, respectively, over the twisted polynomial ring
\[
A[F]=\frac{A\{F\}}{\la a^pF-Fa\mid a\in A \ra},
\]
where $A\{F\}$ is the free non-commutative ring in one variable $F$. 

In either case, the $A[F]$-module structure is encoded by an $A$-module equipped with an additional $A$-linear structure map that corresponds to the respective action of $F$: for Frobenius modules, this is a morphism $\varphi_M:M\to F_*M$, and for Cartier modules, a morphism $C_M:F_*M\to M$. 

Moreover, since the pushforward functor $F_*$ has both a left adjoint $F^*$ and a right adjoint $F^\flat$, these structure maps are equivalent to the corresponding \emph{adjoint structure maps}. These observations are summarized in the following propositions, see \cite[Proposition~2.18]{BBCartierCrstals} and \cite{EmersonKisin}.

\begin{prop}
\label{equiv ch of Cartier}
Let $A$ be an $\mathbb{F}_p$-algebra. A Cartier module over $A$ is given by the following equivalent structures:
\begin{enumerate}
\item A right $A[F]$-module $M$,
\item a pair $(M,C_M)$ consisting of an $A$-module $M$ and an $A$-linear map $C_M:F_*M\to M$, and
\item a pair $(M,\kappa_M)$ consisting of an $A$-module $M$ and an $A$-linear map $\kappa_M:M\to F^{\flat} M$.
\end{enumerate}  
A morphism of Cartier modules is an $A$-linear map that respects these additional structures.
\end{prop}

\begin{prop}
\label{equiv ch of Frobenius}
Let $A$ be an $\mathbb{F}_p$-algebra. A Frobenius module over $A$ is given by the following equivalent structures:
\begin{enumerate}
\item A left $A[F]$-module $M$,
\item a pair $(M,\varphi_M)$ consisting of an $A$-module $M$ and an $A$-linear map $\varphi_M:M\to F_*M$, and
\item a pair $(M,\theta_M)$ consisting of an $A$-module $M$ and an $A$-linear map $\theta_M:F^*M\to M$.
\end{enumerate}  
A morphism of Frobenius modules is an $A$-linear map that respects these additional structures.
\end{prop}

The third perspective offered by the previous propositions leads to the following definition: 

\begin{defn}
A Cartier module (resp. Frobenius module) $M$ is called \emph{unit Cartier module} (resp. \emph{unit Frobenius module}), if its adjoint structure morphism $\kappa_M:M\to F^{\flat}M$ (resp. $\theta_M: F^*M\to M$) is an isomorphism. We denote by $\ucart_A$ (resp. $\ufrob_A$) the full subcategory of unit modules.
\end{defn}

\begin{remark}
In \cite[Definition~2.3, Remark~2.4.]{MaGlobalDimension}, unit Cartier modules are referred to as "unit right $A[F]$-modules". In \cite{LyubeznikFModulesApplicationsToLocalCohomology}, unit Frobenius modules are called "$F$-modules", while in \cite{EmersonKisin} they are called "unit $A[F]$-modules". The terminology "Frobenius module" and "Cartier module" used in this article follows the conventions of \cite{BhattLurie} and \cite{BBCartierCrstals}, respectively.
\end{remark}
As explained above, if $R$ is regular $F$-finite, we have a canonical module $\omega_R$ together with an isomorphism $\kappa:\omega_R\to F^{\flat}\omega_R$. As $\omega_R$ is locally free of rank $1$, tensoring with this module induces an auto equivalence of $\operatorname{Mod}_R$. Moreover, due to the isomorphism $F^\flat(\omega_R \otimes M) \cong F^\flat\omega_R \otimes F^*M \cong \omega_R \otimes F^*M$ (cf. \cite[Corollary~5.8]{BBCartierModules}) one obtains the following equivalence (see also \cite[Theorem~2.5]{MaGlobalDimension}):
\begin{theorem}
\label{unit right is equiv to F-mod}
Let $R$ be a regular $F$-finite ring of positive characteristic $p$. Then an isomorphism $\omega_R \cong F^\flat\omega_R$ yields an equivalence of categories:
\[
\xymatrix{
\ucart_R\ar@/^/[rr]^{\omega_R^{-1}\otimes_R-} & &\ufrob_R \ar@/^/[ll]^{\omega_R\otimes_R-}
}\]
\end{theorem}
As tensoring with an invertible $R$-module affects neither the injective dimension nor the support of the underlying $R$-module, we obtain as an immediate consequence the following analog of \cite[Theorem 1.4]{LyubeznikFModulesApplicationsToLocalCohomology} for unit Cartier modules: 
\begin{corollary}
\label{support inj dim unit cartier}
Let $R$ be a regular $F$-finite ring of positive characteristic $p$. Then 
\[\injdim_{R} N \leq \dim(\Supp_R N)\]
for every unit Cartier module $N$.
\end{corollary}

The category of unit Frobenius modules has good properties when the ring is regular. However, if $R$ is not regular, it fails to be an abelian category due to the non-exactness of $F^*$. In \cite{BhattLurie} the authors introduce a derived category of unit Frobenius modules for any noetherian $F$-finite ring, but it remains unclear in op.cit. if it is the derived category of a suitable abelian category. Quite surprisingly it was observed in \cite{BlickleFink} that unit Cartier modules form an abelian category for any noetherian $F$-finite ring. This makes unit Cartier modules a suitable replacement for unit Frobenius modules beyond the regular case.

\begin{example}
Any ring $A$ of positive characteristic with its usual Frobenius action is a unit Frobenius module. Also all localizations of $A$, and more generally localizations of any unit Frobenius module is again a unit Frobenius module. If $A$ is regular, then local cohomology modules of unit Frobenius modules are unit Frobenius modules, cf. \cite{LyubeznikFModulesApplicationsToLocalCohomology}.

If $A$ is noetherian $F$-finite, then the injective hull of the residue field of any local ring of $A$ carries a unit Cartier module structure: This is due to the fact, that if $E$ is an injective hull of $A/\mathfrak{m}$, then $F^\flat E$ is also an injective hull of $A/\mathfrak{m}$ -- hence it is (non-canonically) isomorphic to $E$. 

However, once we have fixed a quasi-isomorphism $\kappa: \omega^\bullet_A \cong F^!\omega^\bullet_A$ as explained above, we can apply local cohomology to obtain a canonical unit Cartier module structure on $E$:
    \[
        E \simeq R\Gamma_{\mathfrak{m}} \omega^\bullet_A \to[R\Gamma_{\mathfrak{m}}]  R\Gamma_{\mathfrak{m}} F^!\omega^\bullet_A \simeq F^!R\Gamma_{\mathfrak{m}} \omega^\bullet_A \simeq F^!E \cong F^\flat E
    \]
The first and last quasi-isomorphism is local duality, and the middle one is due to the fact that $F^!$ commutes with local cohomology.\footnote{This holds because one can directly verify that $\Gamma_\mathfrak{m} \circ F^\flat = F^\flat \circ \Gamma_\mathfrak{m}$, hence the commutation holds for the derived functors as well.} The final isomorphism follows since $E$ is injective.
\end{example}

\section{Injective dimension of (unit) Cartier modules}
\label{inequalities inj dimension}

In this section we prove our two main (in)equalities and derive from these the proof of Theorem A and Theorem B.
\begin{prop}
\label{inj dim cartier le inj dim module}
Let $A$ be an (arbitrary) $\mathbb{F}_p$-algebra. Then, for any $M\in\cart_A$ we have the following inequality:
\[
\operatorname{inj.dim}_{\cart_A}(M) 
\le 
\operatorname{inj.dim}_{A}(M)+1
\]
\end{prop}
\begin{proof}
Let $M$ be a Cartier module over $A$ with $d=\operatorname{inj.dim}_{A}(M)$. We may assume $d<\infty$ as otherwise, there is nothing to show. As shown in \cite[Proof~of~(a)]{BlickleFink}, for any further Cartier module $N$, there exists a long exact $\Ext$ sequence:
\[
... \to \operatorname{Ext}_{A}^{i}(F_*N,M)\to \operatorname{Ext}_{\cart_A}^{i+1}(N,M) \to \operatorname{Ext}_{A}^{i+1}(N,M)\to ...
\]
This implies that $\operatorname{Ext}_{\cart_A}^{d+2}(N,M)=0$ for all Cartier modules $N$; consequently $\injdim_{\cart_A}M \leq d+1$.
\end{proof}

Now assume that $A$ is a noetherian $F$-finite ring of positive characteristic $p$. It follows from \cite[Lemma~2]{BlickleFink} that the inclusion $\ucart_A\hookrightarrow \cart_A$ admits an exact left adjoint, the \emph{unitalization} $u:\cart_A \to \ucart_A$, which maps a Cartier module $M$ to
\begin{equation}
\label{def unitalization}
    u(M)=\operatorname{colim}(M\to F^{\flat}M \to F^{\flat 2}M \to ...).
\end{equation}

\begin{lemma}
Let $A$ be noetherian $F$-finite, and $M$ be a unit Cartier module over $A$. Then $M$ is injective as an object of $\ucart_A$ if and only if it is injective as an object in the larger category $\cart_A$. 
\end{lemma}
\begin{proof}
Since the unitalization is exact and left adjoint to the inclusion functor, the latter preserves injectives. The reverse implication follows because the inclusion is fully faithful and left exact: 
Suppose $I$ is a unit Cartier module that is injective as a Cartier module, and that $\varphi:M\hookrightarrow N$ is a monomorphism of unit Cartier modules. By assumption, the map
\[
\Hom_{\ucart_A}(N,I)=\Hom_{\cart_A}(N,I)\to \Hom_{\cart_A}(M,I)=\Hom_{\ucart_A}(M,I)
\]
is surjective, proving that $I$ is injective as a unit Cartier module.
\end{proof}

If now $R$ is regular $F$-finite, then Kunz's theorem \cite{KunzCharacterizationsOfRegularLocalRings} shows that $F_*R$ is a projective $R$-module, and consequently, the functor $F^{\flat}$ is exact. A simple application of the snake lemma then shows, that the inclusion $\ucart_R\hookrightarrow \cart_R$ is exact (caution: in general this inclusion is not exact!). With this in mind, we obtain the following lemma:

\begin{lemma}
\label{inj dim unit equals nonunit}
Let $R$ be a regular $F$-finite ring of positive characteristic $p$. For any unit Cartier module $M$ over $R$, we have the following equality of injective dimensions:
\[
\operatorname{inj.dim}_{\ucart_R}(M)=\operatorname{inj.dim}_{\cart_R}(M)
\]
\end{lemma}
\begin{proof}
Since $R$ is regular, the inclusion $\ucart_R\hookrightarrow \cart_R$ is exact. Therefore, as it preserves injectives, it also preserves injective resolutions, proving the inequality $\ge$. Conversely, let $d=\operatorname{inj.dim}_{\cart_R}(M)$, and consider an exact sequence 
\[
0\to M \to I^0 \to I^1 \to ... \to I^{d-1} \to C \to 0
\]
in $\ucart_R$, where $I^0,..., I^{d-1}$ are injective in $\ucart_R$. As argued for the previous inequality, this is also a resolution of $M$ in $\cart_R$ with $I^0,..., I^{d-1}$ remaining injective. By the assumption on the injective dimension, $C$ is an injective object in $\cart_R$. Therefore, the previous lemma implies that it is also injective in $\ucart_R$. Thus, the above sequence is an injective resolution of $M$ in $\ucart_R$.
\end{proof}
Combining Proposition~\ref{inj dim cartier le inj dim module} and Lemma~\ref{inj dim unit equals nonunit} with Corollary~\ref{support inj dim unit cartier} we obtain the following corollary:
\begin{corollary}
\label{injective dimension bounded by support regular case}
Let $R$ be a regular noetherian $F$-finite ring of positive characteristic $p$. For any unit Cartier module $M$ we have the inequality
\[
\operatorname{inj.dim}_{\ucart_R}(M)\le \operatorname{inj.dim}_R(M)+1\le \operatorname{dim}(\operatorname{Supp_R(M)})+1.
\]
\end{corollary}
Using the equivalence of Theorem~\ref{unit right is equiv to F-mod} again, we obtain our main Theorem A:
\begin{theorem}
\label{injdim of F-module}
Let $R$ be a regular $F$-finite ring of prime characteristic $p$. Then 
\[\injdim_{\ufrob_R} (M) \le \operatorname{inj.dim}_R(M)+1\ \leq \dim(\Supp_R M)+1\]
for every unit Frobenius module $M$.
\end{theorem}
The inequality $\operatorname{inj.dim}_R(M)\ \leq \dim(\Supp_R M)$ is the content of Theorem \ref{injdim as R-mod} and it may be the case that $\operatorname{inj.dim}_R(M)\ < \dim(\Supp_R M)$ (for instance, when $M=E_R(R/\mathfrak{p})$ is the injective hull of $R/\mathfrak{p}$ for a non-maximal prime ideal $\mathfrak{p}$). Thus, when $R$ is regular,  $\operatorname{inj.dim}_R(M)+1$ is a tighter bound than $\dim(\Supp_R M)+1$ for $\injdim_{\ufrob_R} (M)$. However, when $R$ is not regular, $\operatorname{inj.dim}_R(M)$ may no longer be finite. Hence it is more feasible to extend the inequality $\injdim_{\ufrob_R} (M) \leq \dim(\Supp_R M)+1$ to singular rings as in the statement of Theorem B.

\section{Reduction to the regular case}
\label{KashiwaraForUnit}
In this section we breifly recall some aspects of the the Kashiwara-type equivalence from \cite{BlickleFink, BBCartierCrstals} for Cartier crystals in the language of unit Cartier modules which is used in the present note. This is then used to extend Corollary~\ref{injective dimension bounded by support regular case} -- which deals with the case that $R$ is regular -- to the general $F$-finite case, concluding the proof of Theorem B. 

Let $\psi: R \to A$ be a map of noetherian $F$-finite rings. The pair of adjoint functors 
\[
\psi_*:\operatorname{Mod}_A
\rightleftarrows 
\operatorname{Mod}_R: \psi^{\flat}
\]
induces adjoint functors on the categories of Cartier modules.
\[
\psi_*:\cart_A \rightleftarrows  \cart_R: \psi^\flat
\]
The unitalization functor $u \colon \cart \to \ucart$ has as its kernel the Serre subcategory of locally nilpotent Cartier modules, i.e. such $M \in \cart$ such that $uM=0$. By \cite[Proposition 6]{BlickleFink}
 unitalization factors through the category $\cryscart$ of Cartier crystals, which is the Serre quotient of the category of Cartier modules by locally nilpotent ones,
\[
\xymatrix{
\cart \ar@{->}[rrd]^{\pi} \ar@{->}[rr]^{u} &  & \ucart \\
 &  & \cryscart \ar@{->}[u]_{\overline{u}}
}
\]
inducing an equivalence of categories $\cryscart \cong \ucart$. As both $\psi_*$ and $\psi^\flat$ preserve local nilpotence it follows from the universal property of localization that these descend to adjoint functors 
\[
\psi^{\unit}_{*}:\ucart_A \rightleftarrows  \ucart_R: \psi^\flat_{\unit}
\]
which satisfy $\psi^{\unit}_* u = u \psi_*$ and $\psi^\flat_{\unit} u = u \psi^\flat$. Note that since $\psi$ is affine, $\psi_*$ is exact and hence so is $\psi^{\unit}_*$. Now the analog of \cite[Proposition 9]{BlickleFink} is the following:
\begin{prop}\label{prop.kashiwaraInjDim}
    Let $\psi \colon R \to A$ be a surjective map of noetherian $F$-finite rings with $R$ regular, then the derived unit of adjunction 
    \[
        \eta \colon \operatorname{id}_{D^+(\cart_R)} \to R\psi^\flat_{\unit} \circ \psi^{\unit}_*
    \]
    is an isomorphism.
\end{prop}
\begin{proof}
    Immediate from \cite[Proposition 9]{BlickleFink} via the equivalence $\cryscart \cong \ucart$.
\end{proof}
\begin{remark}
These adjoint functors can be made explicit as follows. Let $N \simeq I^\bullet$ be an injective resolution in $\ucart_R$. Then by definition of derived functors $R\psi_{\unit }^\flat M \simeq \psi_{\unit}^\flat I^\bullet = \psi_{\unit }^\flat uI^\bullet = u \psi^\flat I^\bullet = \psi^\flat I^\bullet$ where the last equality is due to the easily verified fact that $\psi^\flat$ preserves the unit property. Since the inclusion $\ucart \subseteq \cart$ has an exact left adjoint $u$, each $I^j$ is injective in $\cart_R$. Now, if $R$ is regular, the inclusion $\ucart \subseteq \cart$ is exact, hence $I^\bullet$ is also an injective resolution of the underlying Cartier module of $M$. Hence, at least if $R$ is regular, we have 
\[
R\psi^\flat_{\unit}  M = \psi^!M
\]
such that $R\psi^\flat_{\unit}$ is just $\psi^!$ applied to the underlying Cartier module.\footnote{If $R$ is not regular we still have that $R\psi^\flat_{\unit}M = \psi^! R \operatorname{incl} M$ where $\operatorname{incl} : \ucart_R \to \cart_R$ is the inclusion and $R \operatorname{incl}$ is its right derived functor.}

Similarly we have for $N \in \ucart_A$ that 
\[
\psi^{\unit}_*N = \psi^{\unit}_* uN = u \psi_* N.
\]
Note that $\psi_*$ generally does not preserve the unit property. The adjoint structural map of $\psi_*N$ is given by
\[
\kappa_{\psi_*(M)}:\psi_*(M)\xrightarrow{\psi_*(\kappa_M)} \psi_*F^{\flat}(M)\to F^{\flat}\psi_*(M).
\]
where the second map is induced by the natural base change transformation $\psi_*F^{\flat}\to F^{\flat}\psi_*$. Hence $\psi_*$ preserves unitality if and only if this base change transformation is an isomorphism.
\end{remark}

Now similar as in the proof of (c) in \cite{BlickleFink} we obtain the following lemma:
\begin{lemma}\label{lem.injdimkashiwara}
Let $\psi : R \to A$ be a surjective map with $R$ regular. Then for any unit Cartier module $N$ over $A$ we have
\[
    \injdim_{\ucart_A} N \leq \injdim_{\ucart_R} u\psi_*N
\]
\end{lemma}
\begin{proof}
Let $u\psi_*N = \psi^{\unit}_*N \to I^\bullet$ be an injective resolution of $\psi^{\unit}_*N$ in $\ucart_R$ of length $\leq k$. By definition of right derived functors we obtain a quasi-isomorphism $R\psi^\flat_{\unit}\psi^{\unit}_*N \simeq \psi^\flat_{\unit} I^\bullet$. Composing this with the equivalence $N \simeq R\psi^!_{\unit}\psi^{\unit}_*N$ of Proposition~\ref{prop.kashiwaraInjDim} we get a quasi-isomorphism $N \simeq \psi^\flat_{\unit} I^\bullet$. Since $\psi^\flat_{\unit}$ is right adjoint to the exact functor $\psi_*^{\unit}$, and hence preserves injectivity, we see that $N$ has a injective resolution of lenght $\leq k$ in $\ucart_A$.
\end{proof}
With these preparations at hand we can conclude the proof of Theorem B.
\begin{theorem}
\label{thm: global dim of Cartier crystals2}
Let $A$ be an $F$-finite noetherian ring of prime characteristic $p$. Then 
\[\injdim_{\ucart_A}(N)\leq \dim(\Supp_A N)+1.\]
In particular, the global dimension of the category of unit Cartier modules is at most $\dim(A)+1$.
\end{theorem}
\begin{proof}
Since $A$ is $F$-finite, $A$ admits a surjection $\psi: R\to A$ where $R$ is a regular $F$-finite ring according to \cite[Remark~13.6]{GabberTStructure}. If $N$ is a unit Cartier module over $A$, then via the closed immersion $\Spec A \subseteq \Spec R$ we have $\Supp_A N = \Supp_R \psi_*N$. Since $F^\flat$ commutes with localization we have $\Supp_R \psi_*N \supseteq \Supp_R F^\flat \psi_* N$. The adjoint structural map of $\psi_*N$ is given as follows:
\[
\psi_*N \to[\psi_*\kappa_N] \psi_*F^\flat N \to[b.c.] F^\flat \psi_*N.
\]
The first map is an isomorphism since $\kappa_N$ is an isomorphism as $N$ is unit Cartier. The base change map is given by composition with $F^*\psi$:
\[
\psi_*\Hom_A(F_*A,N) \to[\circ F^*\psi] \Hom_R(F_*R,\psi_*N)
\]
which is injective, since $\psi$ is surjective. Hence $\psi_*N \to F^\flat\psi_*N$ is injective. Combining this with the above we get $\Supp_R \psi_*N = \Supp_R F^\flat \psi_* N$. Iterating this argument we conclude that \[\Supp_A N=\Supp_R \psi_*N = \Supp_R u\psi_*N.\] 

%
Combining the inequality
\[
    \injdim_{\ucart_A} N \leq \injdim_{\ucart_R} u\psi_*N.
\]
from Lemma~\ref{lem.injdimkashiwara} with Corollary~\ref{injective dimension bounded by support regular case} applied to the unit Cartier module $u\psi_*N$ over the regular ring $R$, we get \[\injdim_{\ucart_A} N \leq \dim \Supp_R u\psi_*N + 1=\Supp_A N+1.\] 
\end{proof}

\begin{remark}
Note that $\ucart_A$ is equivalent to the category of quasi-coherent Cartier crystals on $A$ (\cite[Proposition 6]{BlickleFink}). Therefore the same bounds on the injective dimension as in Theorem \ref{thm: global dim of Cartier crystals2} also apply to objects of the category of quasi-coherent Cartier crystals on $A$.
\end{remark}

\begin{remark}
\label{bound is sharp for CM ring}
It is pointed out in \cite[\S3]{BlickleFink} that if $A$ is an equidimensional Cohen-Macaulay ring of prime characteristic $p$ then 
\[\injdim_{\ucart_A}(\omega_A)=\dim(A)+1.\] 
Combining this with Theorem \ref{thm: global dim of Cartier crystals2} shows that: (i) the global injective dimension of $\ucart_A$ is equal to $\dim(A) +1$ in this case; (ii) the bound in Theorem \ref{thm: global dim of Cartier crystals2} is optimal.
\end{remark}

\begin{remark}
    We expect that the regularity assumptions in some of the statements in this section could be removed with some extra care. Since this broader generality is not needed for our goals we leave this as an exercise to the reader.
\end{remark}

\end{document}